# A SURVEY OF CHARACTER SHEAVES

G. Lusztig

**1. Notation.** For a finite group $\Gamma$ we denote by $\mathrm{Irr}\Gamma$ the set of irreducible representations of $\Gamma$ up to isomorphism.

Let $G$ be a connected reductive group over an algebraically closed field $K$. Let $\mathcal{B}$ be the variety of all Borel subgroups of $G$. Let $W$ be the Weyl group of $G$ identified by $w \mapsto \mathcal{O}_w$ with the set of $G$-orbits on $\mathcal{B} \times \mathcal{B}$ (diagonal action). Let $I$ be the set of simple reflections of $W$.

**2.** We now assume that $K$ is an algebraic closure of the finite field $F_q$ and $F : G \to G$ is the Frobenius map corresponding to a split $F_q$-structure on $G$. This induces an endomorphism of $\mathcal{B}$ denoted again by $F$. For any $w \in W$ we set $X_w = \{B \in \mathcal{B}; (B, F(B)) \in \mathcal{O}_w\}$, a smooth subvariety of $\mathcal{B}$ of dimension $l(w)$ (length of $w$) introduced in [DL76]. The subvarieties $X_w$ ($w \in W$) form a partition of $\mathcal{B}$. Each $X_w$ is stable under the conjugation action of the finite group $G^F$. The irreducible representations of $G^F$ which appear in the induced action on the $l$-adic cohomology of $X_w$ are said to be unipotent representations. (The other representations are realized using certain finite coverings of $X_w$.) The unipotent representations of $G^F$ were classified in [L84].

**3.** In a paper [L79] written in 1977, I tried to see what happens if, in the definition of the variety $X_w$, one replaces the Frobenius map by conjugation by $s$, a regular semisimple element of $G$. (Thus, we can define ${}^sY_w = \{B \in \mathcal{B}; (B, sBs^{-1}) \in \mathcal{O}_w\}$.) These varieties again form a partition of $\mathcal{B}$ into pieces which can be shown to be nonempty and smooth. In some sense they can be considered as limits of $X_w$ for $q \mapsto 1$. For example $X_1$ has $q+1$ $F_q$-rational points and ${}^sY_1$ has $1+1$ $F_q$-rational points (at least if $s$ is split). So it can be hoped that ${}^sY_w$ contain some information on $X_w$. (The similarity of ${}^sY_w$, $X_w$ is discussed in [L21].) Now there is no longer an action of an interesting group on ${}^sY_w$ but, if $s$ is allowed to vary, the cohomologies of ${}^sY_w$ form interesting local systems on the variety of regular semisimple elements in $G$. In this paper I showed that from these local systems one can recover information about the unipotent representations $\rho$ of $G^F$ such that $\rho$ appears in the space of functions on the finite set $X_1$. In particular, I could show that for such $\rho$ the multiplicity of $\rho$ in $\sum_i (-1)^i H^i_c(X_w, \bar{\mathbf{Q}}_l)$ is independent of $q$ and in fact is expressible in terms of geometry over $\mathbf{C}$ (since ${}^sY_w$, unlike $X_w$, made sense over $\mathbf{C}$). This paper was for me the beginning of a geometric theory of characters of $G$.

**4.** In 1980, I tried to compute in the case where $G = SL_2(K)$, the intersection cohomology complex $\mathcal{K}$ of $G$ with coefficients in the nontrivial local system of rank 1 on the set of regular semisimple elements given by the natural double covering of this set. To do this, I noticed that this double covering, when extended to the whole of $G$ as the Grothendieck-Springer resolution $\{(g, B) \in G \times \mathcal{B}; g \in B\} \to G$, is a small map (in the sense of Goresky-MacPherson) which allowed me to perform calculations. These calculations showed that the alternating sum of traces of the Frobenius map on the stalks of the cohomology sheaves of $\mathcal{K}$ (the "characteristic function" of $\mathcal{K}$) was giving exactly the values of the character of the $q$-dimensional (or Steinberg) irreducible representation of $SL_2(\mathbf{F}_q)$. This was for me a strong indication that for a general $G$, the characters of irreducible representations of $G^F$ are closely related to the characteristic functions of certain simple perverse sheaves on $G$ and I thus started to look for those perverse sheaves. As a first step, I observed that the Grothendieck-Springer resolution was small for general $G$ and, as an application, I found a new construction of the Springer representations of the Weyl group $W$ in terms of intersection cohomology which, unlike Springer's original definition, made sense in arbitrary characteristic. This appeared in [L81]. At the time (1982) when my book [L84] was written, I had the definition of the required collection of simple perverse sheaves (or character sheaves) in the case where $G$ was $GL_n$, but I knew that for other groups the analogous collection was not complete and I conjectured that it can be completed.

**5.** Later, in 1983, I found two definitions of the collection of character sheaves on a general $G$. One, which appeared in [L84a], was describing the character sheaves as explicit intersection cohomology complexes on certain subvarieties of $G$; this involved a generalization of the Springer representations of $W$ and a generalization of the Grothendieck-Springer resolution which satisfied the definition of a small map except that it was not proper in general (but it could still be used). A second definition was given in a series of papers on character sheaves [L85,L85a,L85b,L86,L86a]. It used the idea of [L79] but applied to not necessarily regular semisimple elements (that series of papers also contains a proof of the equivalence of the two definitions, except for exceptional groups in small characteristic where the proof was completed only in [L12]). The second definition has the virtue that it is in many ways similar to the definition of unipotent representations; in particular the notion of unipotent character sheaf is defined and it is suitable for classification (the classification of character sheaves in the second definition can be essentially reduced to that of unipotent characters sheaves of $G$ and of various smaller groups). On the other hand, the first definition is better suited for explicit calculation of characteristic functions.

**6.** We now give the definition of unipotent character sheaves on $G$ (no assumption on the characteristic of $K$). For each $w \in W$ we define

$$\pi_w : Y_w = \{(g,B) \in G \times \mathcal{B}; (B, gBg^{-1}) \in \mathcal{O}_w\} \to G$$

by $(g,B) \mapsto g$. The simple perverse sheaves on $G$ which appear in the perverse cohomology of $(\pi_w)_!\bar{\mathbf{Q}}_l$ for some $w$ are called the unipotent character sheaves of $G$. (The general character sheaves are obtained by replacing $\bar{\mathbf{Q}}_l$ in $(\pi_w)_!\bar{\mathbf{Q}}_l$ by certain local systems of rank 1 on $Y_w$.)

**7.** Let $Z_G$ be the centre of $G$. Assume that $G/Z_G$ is simple. A unipotent character sheaf $A$ on $G$ is said to be *cuspidal* if the following condition is satisfied: the set of $g \in G$ such that $A|\{g\} \ne 0$ is a union of finitely many conjugacy classes times $Z_G^0$ (it is then automatically a single conjugacy class times $Z_G^0$).

Let $CS(G)$ be the set of (isomorphism classes of) unipotent character sheaves on $G$. Let $CS^\emptyset(G) = \{A \in CS(G); A \text{ cuspidal}\}$.

**8.** Let $\mathrm{Irr}_{spec}(W)$ be the set of special (irreducible) representations of $W$. For any $E \in \mathrm{Irr}_{spec}(W)$ we have defined in [L84] a certain finite group $\Gamma_E$. Let $M(\Gamma_E)$ be the set of all pairs $(x,\sigma)$ where $x \in \Gamma_E$ is defined up to conjugacy and $\sigma \in \mathrm{Irr} Z_{\Gamma_E}(x)$. One of the results in [L85,L85a,L85b,L86,L86a] is the definition of a surjective map

(a) $CS(G) \to \mathrm{Irr}_{spec}(W)$

whose fibre at $E \in \mathrm{Irr}_{spec}(W)$ is in bijection with $M(\Gamma_E)$.

In particular $CS(G)$ can be regarded as a set independent of $K$. Similarly, $CS^\emptyset(G)$ can be regarded as a set independent of $K$.

An interpretation of the map (a) was given in [L15] assuming that the characteristic of $K$ is a good prime or 0.

Let $A \in CS(G)$. There is a unique unipotent class $C$ of $G$ such that:

(a) for any conjugacy class $D$ in $G$ such that $D_u = C$ we have $A|D =$ local system (up to shift);

(b) if $D$ is a conjugacy class in $G$ such that $\dim D_u \ge \dim C$, $D_u \ne C$, then $A|D = 0$.

Here $D_u$ is the set of unipotent parts of elements in $D$. The irreducible representation $E$ of $W$ corresponding to $C$ under Springer correspondence is special. This defines the required map $A \mapsto E$.

**9.** We describe the set $CS^\emptyset(G)$ in the various cases assuming that $G/Z_G$ is simple.

$G$ of type $A_n$, $n \ge 1$: $CS^\emptyset(G) = \emptyset$.

$G$ of type $B_n$ or $C_n$, $n \ge 2$: $CS^\emptyset(G) = \emptyset$ unless $n = k^2 + k$ for some $k \ge 1$ in which case it consists of exactly one object.

$G$ of type $D_n$, $n \ge 4$: $CS^\emptyset(G) = \emptyset$ unless $n = k^2$ for some even $k \ge 2$ in which case it consists of exactly one object.

$G$ of type $E_6$: $CS^\emptyset(G)$ consists of exactly two objects.

$G$ of type $E_7$: $CS^\emptyset(G)$ consists of exactly two objects.

$G$ of type $E_8$: $CS^\emptyset(G)$ consists of exactly 13 objects.

$G$ of type $F_4$: $CS^{\emptyset}(G)$ consists of exactly 7 objects.

$G$ of type $G_2$: $CS^{\emptyset}(G)$ consists of exactly 4 objects.

We now explain the apparition of the numbers $13, 7, 4$ in the last three cases (corresponding to Coxeter diagrams which are "terminal"). In these cases the number of positive roots can be written uniquely as a product of three consecutive integers: $120 = 4.5.6, 24 = 2.3.4, 6 = 1.2.3$ We then have $13 = 4 + 5 + 6 - 2$,7=2+3+4-2, $4 = 1 + 2 + 3 - 2$. Also in these cases $CS^{\emptyset}(G)$ corresponds to a subset of a single $M(\Gamma_E)$ and we have $|\Gamma_E| = 120, 24, 6$ (in fact $\Gamma_E$ is the symmetric group $S_5, S_4$ or $S_3$). Moreover we have $|M(\Gamma_E) = 39, 21, 4$ which are divisible by $13, 7, 4$. (I have no explanation for these facts.)

**10.** For $J \subset I$ let $W_J$ be the subgroup of $W$ generated by $J$ and let $N_W W_J$ be its normalizer. Let $L_J$ be a Levi subgrup of a parabolic subgroup of type $J$.

For $A' \in CS(L_J)$ there is an induced object $ind(A')$ (see [L85,§4]); this is a semisimple perverse sheaf on $G$, (a direct sum of objects in $CS(G)$). If $A' \in CS^{\emptyset}(L_J)$ then the endomorphism algebra of $ind(A')$ can be identified with the group algebra of $N_W W_J/W_J$ hence any $E' \in \text{Irr}(N_W W_J/W_J)$ defines a simple perverse sheaf $A'[E'] = Hom_{N_W W_J/W_J}(E', ind(A')) \in CS(G)$ (a direct summand of $ind(A')$). Let

$$CS'(G) = \{(J, A', E'); J \subset I, A' \in CS^{\emptyset}(L_J), E' \in \text{Irr}(N_W W_J/W_J)\}.$$

Then $(J, A', E') \mapsto A'[E']$ is a bijection $CS'(G) \to CS(G)$.

Let $\mathcal{P}$ be the set of prime numbers. For $r \in \mathcal{P}$ (resp. $r = 0$) let $G^r$ be a connected reductive group over an algebraic closure of the field $F_r$ (resp. over $\mathbf{C}$) of the same type as $G$. For $A \in CS(G)$ let $A^r \in CS(G^r)$ be the object corresponding to $A$ under the bijection $CS(G) \leftrightarrow CS(G^r)$. Assuming that $A' \in CS^{\emptyset}(G)$ (so that $A'^r \in CS^{\emptyset}(G^r)$, let $\mathcal{P}_{A'}$ be the set of all $r \in \mathcal{P}$ such that $\{g \in G^r; A'^r_{\{} g\} \neq 0\}$ is a unipotent class $\gamma_{A'^r}$ of $G^r$ times $Z^0_{G^r}$.

For $r \in \mathcal{P} \sqcup \{0\}$ let $\mathcal{U}_r$ be the set of unipotent classes in $G^r$. It is known that we can regard naturally $\mathcal{U}_0 \subset \mathcal{U}_r$ (Spaltenstein). We can form $\mathcal{U}_* = \cup_r \mathcal{U}_r$ where $\mathcal{U}_0$ is identified with its image in $\mathcal{U}_r$ for any $r \in \mathcal{P}$.

**11.** Following [L23] we define a (surjective) map $CS'(G) \to \mathcal{U}_*$ (hence $CS(G) \to \mathcal{U}_*$) . (We assume that $G/Z_G$ is simple.)

Let $(J, A', E') \in CS'(G)$.

Assume first that $\mathcal{P}_{A'} \neq \emptyset$. For any $r \in \mathcal{P}_{A'}$, the generalized Springer correspndence [L84a] associates to $A'|\gamma_{A'^r}$ and $E'$ a unipotent class $\xi_r \in \mathcal{U}_r$ and a local system on it. We set $t_r(J, A', E') = \xi_r$ viewed as an element of $\mathcal{U}_*$.

There are three cases.

1) $\mathcal{P}_{A'} = \{r_0\}$.

We set $t(J, A', E') = t_{r_0}(J, A', E')$.

2) $\mathcal{P}_{A'} = \mathcal{P}$ and $t_r(J, A', E') \in \mathcal{U}_*$ is independent of $r$.

We set $t(J, A', E') = t_r(J, A', E')$ for any $r \in \mathcal{P}$.

3) $\mathcal{P}_{A'} = \mathcal{P}$ and we are not in case 2. Then there is a unique $r_0 \in \mathcal{P}$ such that $t_r(J, A', E')$ is constant for $r \in \mathcal{P} - \{r_0\}$.

We set $t(J, A', E') = t_{r_0}(J, A', E')$.

Next we assume that $\mathcal{P}_{A'} = \emptyset$. Then $J = I, E = 1$ and $G$ must be of type $E_8$ (there are two possible values for $A'$).

We set $t(J, A', E')$ = regular unipotent class. This completes the definition of $t : CS'(G) \to \mathcal{U}_*$.

**12.** Let $z_G$ be the largest coefficient of the highest root of $G$. (Thus $z_G = 6$ if $G$ has type $E_8$ and $z_G = 4$ if $G$ has type $F_4$.) One can show that the fibre of $t$ at the regular unipotent class is in bijection with $\sqcup_{m \in [1, z_G]} \rho_m$ where $\rho_m$ is the set of primitive $m$-th roots of 1 in $\bar{\mathbf{Q}}_l$.

**13.** We state a conjectural alternative description of the map $t : CS(G) \to \mathcal{U}_*$ assuming that $G/Z_G$ is simple and (for simplicity) that $G$ is of classical type or $E_6, E_7$ or $F_4$ so that $\mathcal{U}_* = \mathcal{U}_2$. We expect that if $A \in CS(G)$ then there is a unique $\gamma \in \mathcal{U}_2$ of maximal dimension such that $A^2|_\gamma \neq 0$ and that $t(A) = \gamma$. (Recall that $A^2 \in CS(G^2)$ corresponds to $A$.) For example, if $A = \bar{\mathbf{Q}}_l$ (up to shift) then $\gamma$ is the regular unipotent class of $G^2$. If $A \in CS(G)$ corresponds to the Steinberg representation then $\gamma = \{1\}$; if $G = Sp_4(K)$ and $A \in CS^{\emptyset}(G)$, then $\gamma$ is the regular unipotent class of $G^2$.

**14.** We now comment on other papers on character sheaves.

In [L87] a notion of character sheaf on a semisimple Lie algebra is introduced and studied (using Deligne-Fourier transform).

In [L86c], [L90],[L92],[L92a],[S95] the connection of character sheaves to characters of irreducible representations of finite reductive groups is studied.

In [L04,L04a,L04b,L04c],[L04d,L04e,L05,L06a,L06b,L09] the notion of character sheaf on a connected component of a possibly disconnected reductive group is introduced and studied.

Let $\mathbf{P}$ be a conjugacy class of parabolic subgroups of $G$. In [L04f, L04g] the notion of "parabolic" character sheaf on the space

$$\{(P, P', gU_P); P \in \mathbf{P}, P' \in \mathbf{P}, g \in G, gPg^{-1} = P'\}$$

is introduced and studied. (An affine analogue is studied in [L10].) When $\mathbf{P} = \{G\}$ one recovers the usual notion of character sheaf.) This is used in [L25] to give a geometric realization of certain Hecke algebras with unequal parameters in terms of perverse sheaves.

In [L06] a conjecture on the existence of a theory of character sheaves on connected unipotent groups is stated. (The conjecture was later proved in [BD13].)

In [L15a] it is shown that the unipotent character sheaves corresponding to a fixed special representation are in natural bijection with the simple objects of the centre of an asymptotic Hecke category (in characteristic zero this was proved earlier in [BFO12] where the case of positive characteristic was conjectured).

In [L15b], [L14] a theory of unipotent character sheaves on loop groups is proposed, based on an affine analogue of the generalized Springer correspondence of [L84a].

In [L17], a theory of (generic) character sheaves on groups of the form $G(\mathbf{k}[\epsilon]/(\epsilon^r))$ is proposed. (It was established there for $r = 2, 4$ and in [BC24] for any $r$.) (This was predicted in [L06].)

Department of Mathematics, M.I.T., Cambridge, MA 02139